\documentclass[]{article}

\usepackage{amsfonts,amsbsy,amscd,amsgen,amsthm,dsfont,amsmath,amssymb,mathrsfs}
\usepackage[colorlinks=true,allcolors=blue]{hyperref}
\usepackage{lipsum}
\usepackage{amsfonts}
\usepackage{graphicx}
\usepackage{epstopdf}
\usepackage{algorithmic}
\usepackage{authblk}

\usepackage{amsopn}

\usepackage{subfigure}
\usepackage{bm} 
\theoremstyle{definition}
\usepackage{xcolor,comment}
\usepackage{enumitem}
\usepackage[utf8]{inputenc}
\usepackage[T1]{fontenc}

\newcommand{\id}{\mathrm d}
\newcommand{\vc}{\mathbf}

\newcommand{\bxi}{\pmb{\xi}}

\renewcommand{\i}{\hat{\mathrm i}}
\renewcommand{\tilde}{\widetilde}
\newcommand{\pard}[2]{\frac{\partial #1}{\partial #2}}

\newcommand{\spn}{\mbox{span}}

\DeclareMathAlphabet\mathbfcal{OMS}{cmsy}{b}{n}

\newcommand{\R}{\mathbb{R}}
\newcommand{\Rn}{\R^n}

\newcommand{\grad}{\nabla} 

\newcommand{\eps}{\varepsilon}

\newtheorem{thm}{Theorem}

\newtheorem{rem}{Remark}
\newtheorem{lem}{Lemma}
\newtheorem{defn}{Definition}
\newtheorem{ass}{Assumption}

\title{Evolution of nonlinear reduced-order solutions for PDEs with conserved quantities}
\author{William Anderson}
\author{Mohammad Farazmand\thanks{Corresponding author's email address: \href{mailto:farazmand@ncsu.edu}{farazmand@ncsu.edu}}}
\affil{Department of Mathematics, North Carolina State University, 2311 Stinson Drive, Raleigh, NC 27695-8205, USA}
\date{}

\begin{document}

\maketitle

\begin{abstract}
Reduced-order models of time-dependent partial differential equations (PDEs) 
where the solution is assumed as a linear combination of prescribed modes 
are rooted in a well-developed theory.
However, more general models where the reduced solutions depend nonlinearly on time varying parameters
have thus far been derived in an ad hoc manner.
Here, we introduce Reduced-order Nonlinear Solutions (RONS): a unified framework for deriving reduced-order models
that depend nonlinearly on a set of time-dependent parameters. 
The set of all possible reduced-order solutions are viewed as a manifold immersed in the function space of the PDE.
The parameters are evolved such that the instantaneous discrepancy between reduced dynamics and the full PDE dynamics is minimized. This results in a set of explicit ordinary differential equations on the tangent bundle of the manifold. In the special case of linear parameter dependence, our reduced equations coincide with the standard Galerkin projection.
Furthermore, any number of conserved quantities of the PDE can readily be enforced in our framework.
Since RONS does not assume an underlying variational formulation for the PDE, it is applicable to a broad class of problems.
We demonstrate the efficacy of RONS on three examples: an advection-diffusion equation, the nonlinear Schr\"odinger equation and Euler's equation for ideal fluids. 
\end{abstract}



\section{Introduction}
\label{sec:Introduction}
Reduced-order models are routinely used to facilitate computational and mathematical analysis of nonlinear partial differential equations (PDEs). The theory is well-developed when the reduced-order solution $\hat u(\vc x,t) = \sum_i q_i(t) u_i(\vc x)$ 
is a linear combination of time-independent modes $u_i$ (see Refs.~\cite{Willcox2015,rowley2017}, for exhaustive reviews). However, a robust mathematical framework is missing when the model $\hat u(\vc x,\vc q(t))$ depends nonlinearly on a set of time-dependent variables $\vc q = (q_1,q_2,$$\cdots,$$q_n)$.

This is despite the fact that such nonlinear reduced-order solutions are ubiquitously used, e.g., 
in prediction of rogue waves~\cite{cousins16,farazmand2017,ruban2015}, 
vortex methods in fluid dynamics~\cite{Beale1985,koumoutsakos_2000}, tracking shocks in supersonic flows~\cite{dalle2010},
and shape optimization~\cite{martins2019,Samareh2001}, to name a few. 
In these studies, the reduced solutions depend nonlinearly on amplitudes, length scales, traveling speed, phase, etc.
Currently, these parameters are evolved in an ad hoc manner based on domain expertise and familiarity with the 
underlying PDE. The purpose of this paper is to propose a unified framework that is broadly applicable to time-dependent nonlinear PDEs.

Nonlinear reduced-order solutions $\hat u(\vc x,\vc q(t))$ are often better suited for efficiently quantifying the dynamics,
as compared to their linear counterpart $\sum_i q_i(t) u_i(\vc x)$.
As a rudimentary example, consider the heat equation in $d$ dimensions,
$\partial_tu= \Delta u$,
and its fundamental solution
\begin{equation}
u(\vc x,\vc q(t)) = A(t)\exp\left(-\frac{|\vc x|^2}{L^2(t)} \right),
\label{eq:heat_fs}
\end{equation}
where $A(t) = (4\pi t)^{-d/2}$ and $L(t) = \sqrt{4t}$.
Here, the parameters are $\vc q = (A,L)$ where the solution depends linearly on the amplitude $A$
and nonlinearly on the length scale $L$. We use this example throughout the paper for illustrative purposes. 
But, more importantly, it already showcases the potency of reduced-order models that depend nonlinearly on time-dependent parameters.
If we were to approximate solution~\eqref{eq:heat_fs} as the linear combination of time-independent modes $u_i$, i.e., 
$\sum_i q_i(t) u_i(\vc x)$, several modes would be required to obtain a reasonable approximation.
On the other hand, the solution can be expressed with a single mode~\eqref{eq:heat_fs} when allowing nonlinear dependence on time-dependent parameters.

Here, we develop a unified framework for evolving the reduced-order solutions $\hat u(\vc x,\vc q(t))$ which depend nonlinearly on time-dependent parameters $\vc q(t)$. For a prescribed dependence on the spatial variable $\vc x$, we view $\hat u$ as a map from the parameter space to the function space $H$ of the PDE. The image of this map is a manifold immersed in $H$. We require the evolution of 
the parameters $\vc q(t)$ to minimize the instantaneous discrepancy between the reduced-order dynamics and full dynamics of the PDE.
This leads to a set of nonlinear ordinary differential equations (ODEs) for $\vc q(t)$, defined on the tangent bundle of the manifold.
In our framework, any number of conserved quantities of the PDE can be readily enforced. Since it does not rely on the variational
formulation of the PDE, our method is broadly applicable. We also show that when $\hat u$ depends linearly on the parameter $\vc q$, 
our reduced-order equations coincide with the standard Galerkin projection models.
We refer to the proposed method as RONS: Reduced-Order Nonlinear Solutions.

We point out that the method presented here should not be confused with \emph{parametric model reduction}~\cite{Willcox2015} where the PDE itself depends on some parameters. In contrast, in RONS, the parametric dependence appears in the reduced-order solution.
Furthermore, the term \emph{nonlinear model reduction} is often used for reduced-order modeling of nonlinear PDEs even when
the reduced solution is a linear combination of modes. In this paper, in addition to the PDE being nonlinear, the 
reduced solution is also a nonlinear function of the reduced variables.

\subsection{Related work}
This paper introduces a novel and unified framework for evolving time-dependent parameters of reduced-order nonlinear solutions. However, as mentioned earlier, 
the idea of using such nonlinear ansatz is not new. Here, we review some of the work which use nonlinear ansatz and the respective ad hoc methods for
evolving their parameters.

The propagation of optical beams in nonlinear dispersive media has been studied by considering the evolution of a 
localized Gaussian ansatz with time-dependent amplitude, length scale, and phase~\cite{berge1994,Desaix91,PerezGarcia1996}. A similar approach has also been taken in the context of nonlinear water waves by 
Ruban~\cite{ruban2015,ruban2015b} and Adcock et al.~\cite{Adcock12,adcock09}. 
Following~\cite{Desaix91}, Ruban~\cite{ruban2015,ruban2015b} uses the Lagrangian associated with the nonlinear Schr\"odinger equation (NLSE) to derive reduced-order equations
for the parameters of the ansatz. Adcock et al.~\cite{Adcock12,adcock09}, on the other hand, use the conserved quantities of NLSE to evolve the parameters.
Cousins and Sapsis~\cite{cousins15} consider the modified NLSE which does not have a known Lagrangian structure or as many conserved quantities as NLSE.
They use a hyperbolic secant as the ansatz, take a second time derivative of the modified NLSE, and project the resulting equation on the subspace of the ansatz to derive their reduced-order equations.

In fluid dynamics, vortex methods decompose the flow field into a combination of vortices, with a prescribed smooth profile, whose centers evolve over time~\cite{chorin_1973,Beale1985}. The motion of each vortex center is determined by computing the induced velocity by the other vortices.

Self-similar solutions of PDE are also obtained by varying the parameters of a prescribed ansatz~\cite{barenblatt_1996}.
Rowley et al.~\cite{Rowley_2003} use symmetry reduction from geometric mechanics to derive the symmetry reduced governing equations. 
Although rigorous, this method is only applicable for reducing continuous symmetries as opposed to RONS which is a more general reduced-order modeling framework.

Finally, we point out that the idea of optimally time-dependent modes (OTD) has been introduced in the context of stability analysis~\cite{otd}. However, OTD is only applicable to linearized PDEs. In addition, to compute each OTD mode, an auxiliary PDE must be solved which renders this method computationally expensive~\cite{babaee17,PRE2016}.
In contrast, RONS is applicable to fully nonlinear PDEs and reduces the computational cost as it only requires solving a relatively small set of ODEs.

\subsection{Outline of the paper}
This paper is organized as follows. In section \ref{sec:setup}, we introduce the problem set-up in which our method is applicable. Section \ref{sec:varmodeling} describes RONS in detail and contains our main results. In section \ref{sec:examples}, we present numerical results on three different examples. Section \ref{sec:conclusion} contains our concluding remarks.

\section{Set-up and preliminaries}
\label{sec:setup}

We consider PDEs of the form
\begin{equation}
\frac {\partial u} { \partial t } =  F(u), \quad u( \vc x, 0) = u_0( \vc x), 
\label{eq:General_PDE}
\end{equation}
for the map $u:D\times \R^+\to \mathbb R^p, (\vc x,t) \mapsto u(\vc x,t)$
where $D$ is a subset of $\mathbb R^k$. Here, $F$ is a potentially nonlinear differential operator.
For simplicity, we restrict our attention to the case $p=1$, where $u$ is a scalar function of the spatial variable $\vc x\in D$ and time $t\geq 0$. Our results generalize to the case $p>1$ in a straightforward fashion. We assume that for any time $t$, $u(\cdot,t)$ belongs to a Hilbert function space $H$ with the inner product $\langle \cdot,\cdot \rangle_H$
and the corresponding norm $\|\cdot \|_H$. The appropriate boundary conditions for the PDE are encoded in the Hilbert space $H$.

We seek approximate solutions of the form $\hat u(\vc x,\vc q(t))$ to PDE~\eqref{eq:General_PDE}. 
The dependence of $\hat u$ on the spatial variable $\vc x$ is prescribed based on the type of the PDE and requires familiarity with its solutions. 
As such, we refer to the approximation solution $\hat u(\vc x,\vc q(t))$ as the \emph{ansatz}.
For instance, the heat kernel discussed in Section~\ref{sec:Introduction} has a Gaussian shape with time-dependent amplitude $A(t)$ and length scale $L(t)$ so that $\vc q(t) = (A(t),L(t))$.
Given the prescribed spatial shape of the ansatz $\hat u$, our objective is to determine the evolution of the parameters $\vc q(t)$, so that the ansatz $\hat u(\vc x,\vc q(t))$
best approximates a true solution $u(\vc x,t)$ of the PDE. Since the ansatz solutions are typically smooth in the spatial variable $\vc x$, we only consider strong solutions of the PDE, assuming their existence and uniqueness.

\begin{figure}
	\centering
	\includegraphics[width=\textwidth]{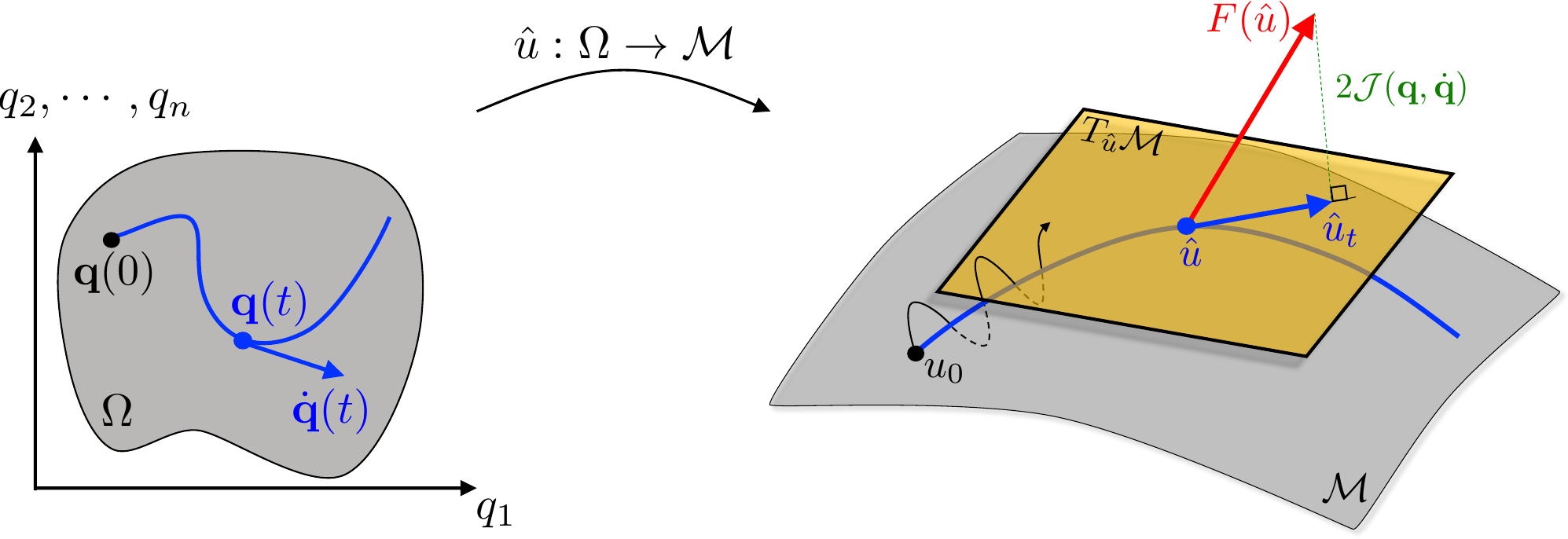}
	\caption{Geometric illustration of RONS. The reduced solution evolves on the ansatz manifold $\mathcal M\subset H$, which is the image of the ansatz map $\hat u$. The blue curve illustrates the evolution of the reduced-order solution. The black curve marks the true trajectory of the PDE which can leave $\mathcal M$ since the manifold is not invariant under the full dynamics.}
	\label{fig:Geometric_Explanation}
\end{figure}

We view the ansatz $\hat u(\vc x,\vc q)$ as a map from the parameters $\vc q\in \Omega\subseteq \Rn$ to the function space $H$,
\begin{align}
\hat u : &\ \Omega \to H\nonumber\\
&\ \vc q \mapsto \hat u(\cdot,\vc q),
\label{eq:ansatz_map}
\end{align}
where $\Omega$ is a simply connected open subset of $\Rn$ over which the ansatz $\hat u$ is well-defined.
For instance, for the heat kernel discussed in Section~\ref{sec:Introduction}, $\Omega$ is the positive quadrant of $\R^2$ 
where $A>0$ and $L>0$. We note that the domain of the map $\hat u$ is over the finite-dimensional parameters $\vc q = (q_1,q_2,\cdots, q_n)$, 
whereas, for a given parameter $\vc q$, its image is a function in the Hilbert space $H$ (see figure~\ref{fig:Geometric_Explanation} for an illustration). We refer to the image of the map $\hat u$ as the ansatz manifold, motivating the following definition.

\begin{defn}[Ansatz Manifold]
	We refer to the set 
	\begin{equation}
	\mathcal M := \left\{u\in H: \exists\, \vc q\in\Omega\; \mbox{with}\; u = \hat u (\cdot,\vc q) \right\}\subset H,
	\end{equation}
	which is the image of the map~\eqref{eq:ansatz_map}, as the ansatz manifold.
\end{defn}

For well-posedness of the reduced-order equations governing the evolution of the parameters $\vc q$, a number of mild assumptions 
regarding the ansatz map $\hat u$ are required. We list these assumptions below and explain their relevance.
\begin{ass}
	We assume that the map $\hat u: \Omega \to H$ has the following properties.
	\begin{enumerate}[label=(H\arabic*)]
		\item It is injective and at least once continuously differentiable with respect to the parameters $\vc q$.
		\item For every parameter $\vc q\in\Omega$, the ansatz $\hat u(\vc x,\vc q)$ is $\ell$-times continuously differentiable with respect to the spatial variable $\vc x$.
		Here, $\ell\in\mathbb N$ is the highest order of spatial derivatives appearing in the PDE~\eqref{eq:General_PDE}.
		\item The map $\hat u: \Omega \to H$ is an immersion.
	\end{enumerate}
	\label{ass:ansatz_map}
\end{ass}

Our reduced-order equations involve partial derivatives of the ansatz with respect to the parameters. Assumption (H1) ensures that these derivatives exist.
Assumption (H2) ensures that the ansatz would in fact approximate strong solutions of the PDE.
Finally, assumption (H3) requires the map $\hat u$ to be an immersion, which essentially means that the tangent space of the ansatz manifold at every base point $\hat u\in \mathcal M$ 
is full rank. More precisely, the tangent space $T_{\hat u}\mathcal M$ at a point $\hat u\in\mathcal M$ is the linear vector space
\begin{equation}
T_{\hat u}\mathcal M := \spn \left\{\pard{\hat u}{q_1}, \pard{\hat u}{q_2},\cdots,\pard{\hat u}{q_n} \right\}.
\end{equation}

The map $\hat u$ being an immersion implies that the partial derivatives $\partial \hat u/\partial q_i$ are linearly independent 
functions and therefore $\dim(T_{\hat u}\mathcal M)=n$ for all $\hat u\in\mathcal M$.
As such, the set $\mathcal M$ is an $n$-dimensional manifold immersed in the function space $H$.
Most importantly, assumption (H3) ensures that the metric tensor, to be defined in Section~\ref{sec:varmodeling}, is invertible.
It should be clear that, as a special case, the ansatz map includes linear superposition of potentially parameter-dependent modes,
\begin{equation}
\hat u (\vc x,\vc q) = \sum_i a_i u_i(\vc x,\pmb\alpha_i),
\end{equation}
where $\vc q$ consists of the parameters $a_i$ and $\pmb\alpha_i$.

As reviewed in the Introduction, earlier studies which use nonlinear reduced-order solutions are confronted with two main tasks,
\begin{enumerate}
	\item Choosing an ansatz manifold, i.e., choosing the shape of the ansatz $\hat u(\vc x,\vc q(t))$,
	\item Choosing a strategy for evolving the variables $\vc q(t)$ in time.
\end{enumerate}
The first task relies on domain expertise and familiarity with the solutions of the PDE. The second task has thus far been carried out based on ad hoc methods.
The present paper introduces a rigorous and unified approach for addressing the second task. Choosing an appropriate ansatz manifold remains a challenging problem
and still relies on domain expertise.

\section{Evolution of reduced-order solutions}
\label{sec:varmodeling}
This section contains our main results. First, in Section~\ref{sec:rom_unconst}, we discuss the minimization principle from which the reduced-order equations are derived. For clarity, we do not assume any conserved quantities for the PDE in this section. 
Next, we show that our reduced-order equations are equivalent to the standard Galerkin projection if the ansatz is a  linear function of the parameters
(Section~\ref{sec:galerkin_equiv}). Finally, for PDEs with conserved quantities, 
in Section~\ref{sec:rom_const} we derived reduced-order equations that respect those conservation laws.

\subsection{Reduced-order equations without conserved quantities}\label{sec:rom_unconst}
Given an ansatz $\hat u(\vc x,\vc q(t))$, we would like to determine an evolution equation for the parameters $\vc q(t)$ such
that the ansatz closely approximates a true solution $u(\vc x,t)$ of the PDE~\eqref{eq:General_PDE}. Since the true solutions are a priori unknown, this cannot be achieved by minimizing some distance metric between the ansatz and the true solution.
Therefore, an alternative metric needs to be used. 

Here, we evolve the ansatz such that its instantaneous dynamics best match the true dynamics of the PDE.
To this end, we consider the \emph{instantaneous error},
\begin{equation}
\mathcal J( \vc q, \dot{ \vc q }) = \frac{1}{2} \|\hat u_t-F(\hat u)\|_H^2,
\label{eq:costfunctional}
\end{equation}
where $\hat u_t$ is shorthand for the derivative of the ansatz $\hat u(\vc x,\vc q(t))$ with respect to time, i.e., 
\begin{equation}
\hat u_t (\vc x,\vc q(t))= \pard{\hat u}{q_i}(\vc x,\vc q(t))\dot q_i(t),
\label{eq:ut}
\end{equation}
where summation over repeated indices is implied. We note that $\hat u_t$ belongs to the tangent space $T_{\hat u}\mathcal M$
of the ansatz manifold $\mathcal M$ at the point $\hat u$.

The cost function~\eqref{eq:costfunctional} measures the instantaneous discrepancy between the dynamics of the ansatz, i.e. $\hat u_t$,
and the dynamics dictated by the PDE, i.e. $F(\hat u)$. Given a set of parameters $\vc q(t)$, at any time $t\geq 0$, 
we seek $\dot{\vc q}(t)$ such that the instantaneous error $\mathcal J(\vc q,\dot{\vc q})$ is minimized.
Figure~\ref{fig:Geometric_Explanation} shows a geometric illustration of the instantaneous error $\mathcal J$.

\begin{rem}
	As an alternative to the instantaneous error, one might be tempted to minimize the \emph{finite-time error},
	\begin{equation}
	S[ \vc q]= \int_{0}^{T} \mathcal J(\vc q(t),\dot{\vc q}(t)) \id t,
	\label{eq:LagrangianEL}
	\end{equation}
	over a time interval $t\in[0,T]$.
	At first sight, this functional seems more suitable, compared to its instantaneous counterpart~\eqref{eq:costfunctional}, since
	it measures the accumulated error between the true dynamics and the ansatz dynamics. 
	However, as we show in Appendix~\ref{sec:EL_unstable}, it generally leads to unstable reduced-order models for the evolution of 
	the parameters. Therefore, in this paper we choose to minimize the instantaneous error.
	\label{rem:FTerror}
\end{rem}

As mentioned earlier, given the parameter $\vc q(t)$ at time $t$, we determine its dynamics by minimizing the instantaneous error~\eqref{eq:costfunctional} over all possible $\dot{\vc q}(t)$. 
More precisely, we solve the minimization problem
\begin{equation}
\min_{\dot{\vc q} \in \R^n} \mathcal J (\vc q, \dot{\vc q}).
\label{eq:gradJ0}
\end{equation} 
We prove that this problem has a unique minimizer which satisfies a first-order ODE for the parameters $\vc q$.
By solving this ODE, the optimal evolution of the parameters $\vc q(t)$ can be obtained.
To prove these results, we need the following lemma.

\begin{lem}[Metric Tensor]
	Let Assumption~\ref{ass:ansatz_map} hold. Then the metric tensor  $M$ defined by
	\begin{equation}
	M_{ij} = \left\langle \pard{\hat u}{q_i}, \pard{\hat u}{q_j}\right\rangle_H,\quad i,j\in\{1,2,\cdots,n\},
	\label{eq:M}
	\end{equation}
	is a symmetric positive-definite matrix for all $\vc q\in \Omega$.
	\label{lem:MSPD}
\end{lem}
\begin{proof} It is clear from its definition that $M$ is symmetric. We show that the matrix is positive definite by proving that 
	$\langle \bxi,M\bxi\rangle>0$ for all nonzero $\bxi\in\Rn$, where $\langle\cdot,\cdot\rangle$ denotes the Euclidean inner product. 
	We first note that
	\begin{equation}
	\langle \bxi,M\bxi\rangle = \left\langle \pard{\hat u}{q_i}\xi_i, \pard{\hat u}{q_j}\xi_j\right\rangle_H  = \left\|\pard{\hat u}{q_i}\xi_i\right\|_H^2\geq 0,
	\label{eq:MatrixM}
	\end{equation}
	where summation over repeated indices is implied.
	Next we show that $\langle \bxi,M\bxi\rangle$ is in fact strictly positive.
	
	Assume that there exists $\bxi\neq\vc 0$ such that $\langle \bxi,M\bxi\rangle=0$. This, together with equation~\eqref{eq:MatrixM}, implies that
	$$\pard{\hat u}{q_i}\xi_i=0.$$
	Since $\bxi\neq \vc 0$, this implies that $\partial\hat u/\partial q_i$, $i=1,2,\cdots, n$, are linearly dependent. However, this violate assumption (H3) that
	the map $\hat u$ is an immersion. Therefore, we must have $\langle \bxi,M\bxi\rangle>0$ which completes the proof.
\end{proof}

The fact that the metric tensor $M$ is symmetric positive-definite, and hence invertible, plays an important role in deriving our reduced-order equations.
Now we state the main result of this section which establishes that there exists a unique solution to the minimization problem \eqref{eq:gradJ0}, and that the minimizer satisfies an explicit ODE.
\begin{thm}\label{thm:qdot_unconst}
	Let Assumption~\ref{ass:ansatz_map} hold. Then there exists a unique solution to the minimization problem
	$\eqref{eq:gradJ0}$. Furthermore, the minimizer satisfies
	\begin{equation}
	\dot{ \vc q } = M^{-1}(\vc q) \vc f(\vc q)
	\label{eq:qdot_unconst}
	\end{equation}
	where $M$ is the metric tensor defined in Lemma \ref{lem:MSPD}
	and $\vc f:\Rn\to \Rn$ is a vector field defined by
	$$f_i = \left\langle \pard{\hat u}{q_i}, F(\hat u)\right\rangle_H,\quad i=1,2,\cdots, n.$$
\end{thm}

\begin{proof}
	First note that the instantaneous error \eqref{eq:costfunctional} can be written more explicitly as
	\begin{align*}
	\mathcal J (\vc q, \dot{ \vc q })  & = \frac{1}{2}\left[ \langle \hat u_t,\hat u_t\rangle_H -2 \langle \hat u_t,F(\hat u)\rangle_H + \langle F(\hat u),F(\hat u)\rangle_H \right] \nonumber\\
	& = \frac{1}{2} \dot q_i\dot q_j \left\langle \pard{\hat u}{q_i}, \pard{\hat u}{q_j}\right\rangle_H - \dot q_i \left\langle \pard{\hat u}{q_i}, F(\hat u)\right\rangle_H + \frac{1}{2} \langle F(\hat u),F(\hat u)\rangle_H\nonumber\\
	& = \frac{1}{2} \langle \dot{\vc q},M(\vc q) \dot{\vc q}\rangle - \langle \dot{\vc q}, \vc f(\vc q)\rangle + \frac{1}{2} \langle F(\hat u),F(\hat u)\rangle_H.
	\end{align*}
	Therefore, $\mathcal J$ is a smooth quadratic function of $\dot{\vc q}$. Furthermore, since the metric tensor $M(\vc q)$ is symmetric positive-definite, 
	$\mathcal J$ is a strictly convex function of $\dot{\vc q}$. As such, it has a unique minimizer satisfying $\grad_{\dot{ \vc q }} \mathcal J (\vc q, \dot{ \vc q })=0$, or equivalently
	$M(\vc q) \dot{ \vc q } - \vc f(\vc q) = 0$.
	Since the metric tensor $M(\vc q)$ is invertible, we obtain equation~\eqref{eq:qdot_unconst} for $\dot{\vc q}$.
\end{proof}

Figure~\ref{fig:Geometric_Explanation} depicts the geometric interpretation of Theorem~\ref{thm:qdot_unconst}. 
Consider an arbitrary, but smooth, evolution of the parameters $\vc q(t)$. This can be viewed as 
a curve $\vc q:[0,\infty)\to \Omega,\, t\mapsto \vc q(t)$ parametrized by time $t$. 
The tangent vector to the curve at a point $\vc q(t)$ is given by $\dot{\vc q}(t)$.
The ansatz $\hat u$ maps this curve onto a curve on the ansatz manifold $\mathcal M$ in the function space $H$.
The tangent vector to this second curve is $\hat u_t$ (see equation~\eqref{eq:ut}) which is bound to belong to the 
tangent space $T_{\hat u}\mathcal M$ of the ansatz manifold, i.e., $\hat u_t \in T_{\hat u}\mathcal M$.

On the other hand, $F(\hat u)$ does not necessarily belong to $T_{\hat u}\mathcal M$ since the ansatz manifold $\mathcal M$ is not 
invariant under the dynamics of the PDE~\eqref{eq:General_PDE}. Therefore, typically there is no evolution of the parameters $\vc q(t)$
such that the ansatz $\hat u(\cdot,\vc q(t))$ would solve the PDE.  If the parameters $\vc q(t)$ evolve according to equation~\eqref{eq:qdot_unconst}, Theorem~\ref{thm:qdot_unconst} ensures that the deviation of the ansatz dynamics $\hat u_t$ from the true 
dynamics $F(\hat u)$ is instantaneously minimized. 
In fact, $\hat u_t$ is the orthogonal projection of $F(\hat u)$ onto the tangent space of the manifold at point $\hat u$.

We close this section with a remark about the initial condition of the parameters.
In order to numerically integrate the reduced-order equation~\eqref{eq:qdot_unconst}, we need to supply an initial condition $\vc q(0) = \vc q_0$.
If the initial condition $u_0$ of the PDE belongs to the ansatz manifold $\mathcal M$, there exists a unique $\vc q_0\in\Omega$ such that 
$u_0(\vc x) = \hat u(\vc x,\vc q_0)$. In this case, there is no ambiguity regarding the appropriate initial parameter values $\vc q_0$.
However, if $u_0\notin\mathcal M$, a criterion must be devised to determine a suitable initial condition $\vc q_0$.
One such criterion, for instance, is to solve the optimization problem,
\begin{equation}
\vc q_0 = \arg\min_{\vc q\in \Omega} \| u_0 - \hat u(\cdot,\vc q)\|_H,
\label{eq:u0}
\end{equation}
which returns the closest point on the ansatz manifold to the initial condition $u_0$.
Optimization problem~\eqref{eq:u0} can potentially be non-convex and computationally expensive to solve, 
but it is solved only once at the initial time.

\subsection{Relation to Galerkin projection}\label{sec:galerkin_equiv}
\label{sec:Galerkin}
In this section, we show that if the ansatz $\hat u(\vc x,\vc q)$ is linear in the parameters $\vc q$, then the reduced-order equations~\eqref{eq:qdot_unconst}
coincides with the standard Galerkin truncation. 

Consider the ansatz
\begin{equation}
\hat u(\vc x,\vc q (t)) = \sum_{i = 1}^n q_i(t) u_i(\vc x),
\label{eq:linAnsatz}
\end{equation}
which is a linear superposition of the prescribed \emph{modes} $\{u_i\}_{i=1}^n$. Without loss of generality, we assume that these modes
are orthonormal with respect to the inner product on the Hilbert space $H$, i.e.,
$\langle u_i,u_j\rangle_H = \delta_{ij}, \quad i,j\in\{1,2,\cdots,n\}$,
where $\delta_{ij}$ denotes the Kronecker delta. These modes span an $n$-dimensional linear subspace $V$ of the function space $H$, where
$V:=\spn\{u_i\}_{i=1}^n$. The subspace $V$ is not generally an invariant subspace under the dynamics of the PDE~\eqref{eq:General_PDE}. 
Therefore, the Galerkin ansatz~\eqref{eq:linAnsatz} may not be an exact solution of the PDE. More precisely, although $\hat u_t$ belongs to $V$,
the right-hand side of the PDE $F(\hat u)$ does not generally lie in $V$. 

To obtain an approximate solution, one defines the projection operator $P:H\to V$ which is an orthogonal projection onto the subspace $V$, and replaces the right-hand side
with $PF(\hat u)$. Substituting the ansatz in the truncated PDE, $\partial_t \hat u = PF(\hat u)$, and taking the inner product with a mode $u_k$, we finally obtain the standard Galerkin projection,
\begin{equation}
\dot q_k = \langle u_k,F(\hat u)\rangle_{H},\quad k=1,2,\cdots, n.
\label{eq:galerkin_linPDE}
\end{equation}
\begin{figure}
	\centering
	\includegraphics[width=0.6\textwidth]{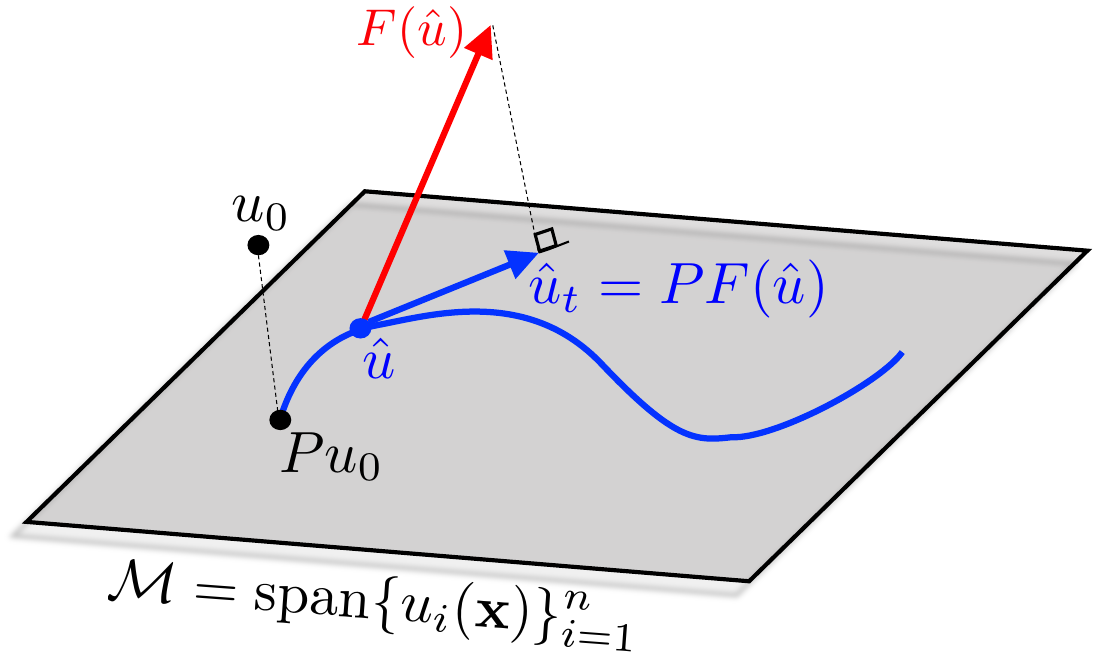}
	\caption{The geometry of the reduced-order equation for a linear ansatz. In this case, the ansatz manifold $\mathcal M$ is a linear subspace of the 
		function space $H$. The reduced-order equations coincide with the Galerkin truncation of the PDE to the modes $u_i$ that span this subspace.}
	\label{fig:galerkin}
\end{figure}

Now, we show that the reduced-order equation~\eqref{eq:qdot_unconst} of Theorem~\ref{thm:qdot_unconst} coincides with the Galerkin truncation~\eqref{eq:galerkin_linPDE}.
First note that, for the linear ansatz~\eqref{eq:linAnsatz}, we have
\begin{equation}
\pard{\hat u}{q_i} = u_i,
\end{equation}
and therefore the metric tensor $M$ is an identity matrix since $M_{ij} = \langle u_i, u_j\rangle=\delta_{ij}$.
Furthermore, for a linear ansatz, the vector field $\vc f:\Rn\to \Rn$ is given by 
\begin{equation}
f_i = \langle u_i,F(\hat u)\rangle_H.
\end{equation}
Substituting these in the reduced-order equation~\eqref{eq:qdot_unconst}, we see that it coincides exactly with the Galerkin truncation~\eqref{eq:galerkin_linPDE}.

Figure~\ref{fig:galerkin} depicts the geometric meaning of the Galerkin projection. In the case of a linear ansatz, the ansatz manifold $\mathcal M$ becomes a linear subspace of 
the function space $H$. The Galerkin method evolves the ansatz so that instantaneously the tangent vector to the path of the ansatz is the orthogonal projection of $F(\hat u)$ 
onto the ansatz subspace. Therefore, our method is a generalization of the Galerkin method to the case when the ansatz manifold is a nonlinear subset of the function space. In this more
general case, the evolution of the ansatz is obtained by projecting $F(\hat u)$ onto the tangent space of the manifold $\mathcal M$ at the ansatz $\hat u$ (see figure~\ref{fig:Geometric_Explanation}).

\subsection{Reduced-order equations with conserved quantities}\label{sec:rom_const}
\label{sec:conserved_quants}
Certain PDEs possess conserved quantities which are invariant along their trajectories~\cite{carlberg2015}. 
Sometimes these conserved quantities are evident because the PDE is derived from underlying conservation laws~\cite{lax1973}.
Other conserved quantities, such as helicity in Euler's equation for ideal fluids~\cite{moffatt1969}, are not evident and are only discovered by further mathematical analysis
of the PDE. In either case, constructing reduced-order models that preserve these conserved quantities is attractive for two main reasons.
First, conserved quantities reflect symmetries and physical properties of the system that reduced-order models should also  
exhibit~\cite{holmes2012}. Secondly, it is known that reduced-order models that violate the conservation laws can develop spurious finite-time blowups
and other non-physical dynamics~\cite{majda2012}. 

In this section, we modify the reduced-order equations developed in Section~\ref{sec:setup} to enforce conserved quantities of the PDE.
Let $I_k: H \to \R$ with $k\in\{1,2,\cdots,m\}$ denote $m$ conserved quantities of the PDE~\eqref{eq:General_PDE}. More precisely, if $u(\cdot, t)$
denotes a solution of the PDE, then $I_k(u(\cdot,t)) = I_k(u_0)$ and therefore $I_k$ is independent of time.
We seek reduced-order equations for evolving an ansatz $\hat u(\vc x, \vc q(t))$ such that the conserved quantities $I_k$ are also invariant 
along the trajectories of the ansatz. For notational simplicity, we write $I_k(\vc q)$ instead of $I_k(\hat u(\cdot,\vc q))$ and view $I_k(\vc q)$
as a map from the parameter space $\Omega\subseteq \Rn$ to the real line.

\begin{rem}
	To be precise, $I_k(\hat u(\cdot,\vc q))$ is the pullback of $I_k:H\to \mathbb \R$ under the ansatz map $\hat u:\Omega \to H$.
	In other words, we should write 
	$$\hat u^\ast I_k(\vc q) = I_k(\hat u(\cdot,\vc q)), \quad \forall \vc q\in \Omega,$$
	where $\hat u^\ast I_k$ denotes the pullback of $I_k$ under the ansatz map $\hat u$. 
	However, this unnecessarily complicates the notation. Therefore, we simply write $I_k(\vc q)$ instead of $\hat u^\ast I_k(\vc q)$.
\end{rem}

To enforce the conserved quantities, we add them as constraints to the optimization problem~\eqref{eq:gradJ0}, and solve
\begin{align}\label{eq:constJ}
& \min_{\dot{\vc q} \in \R^n} \mathcal J(\vc q,\dot{\vc q}),\nonumber\\
&\mbox{s.t.}\quad I_{ k }(\vc q (t))= I_{ k,0 } , \quad k =1,2,...,m, \quad \forall t\geq 0,
\end{align}
where $I_{k,0}$ are prescribed constants. In practice, these constants are determined by the initial condition so that $I_{k,0} = I_k(\vc q_0)$.
To avoid degenerate solutions of the optimization problem, we must make the following assumption.
\begin{ass}
	We assume that, for each $\vc q\in \Omega$, the gradients $\nabla I_1(\vc q), \nabla I_2(\vc q)$, $\cdots$, and $\nabla I_m(\vc q)$ are linearly independent.
	\label{ass:gradI}
\end{ass}

In order to derive the reduced-order equations corresponding to the constrained optimization problem~\eqref{eq:constJ}, we need the following lemma.
\begin{lem}[Constraint Matrix]
	Let Assumptions~\ref{ass:ansatz_map} and~\ref{ass:gradI} hold and
	$M\in \R^{n\times n}$ denote the metric tensor defined in Lemma \ref{lem:MSPD}.  
	Then the \emph{constraint matrix} $C(\vc q)\in \R^{m\times m}$ with entries defined by 
	\begin{equation}
	C_{ij} = \langle \grad I_j ,M^{-1} \grad I_i\rangle,\quad i,j\in\{1,2,\cdots,m\}, 
	\label{eq:constMat}
	\end{equation}
	is symmetric positive-definite for all $\vc q\in \Omega$.
	\label{lem:ASPD}
\end{lem}

\begin{proof}
	First, we define the $n\times m$ matrix $B := \left( \grad I_1|  \grad I_2 |\cdots |\grad I_m \right)$ and note that the constraint matrix can be written as $C = B^\top M^{-1}B$.
	Since the matrix tensor $M$ is symmetric, so is the constraint matrix $C$.
	
	Now we prove that $C$ is also positive-definite. For any nonzero $\vc v\in\R^m$, we have
	\begin{equation}
	\langle  \vc v, C \vc v \rangle = \langle  \vc v, B^\top M^{-1}B  \vc v \rangle  =  \langle B \vc v, M^{-1}B  \vc v \rangle  > 0,
	\end{equation}
	where the last inequality is justified by the facts that $M^{-1}$ is symmetric positive-definite and $B\vc v$ is nonzero for $\vc v \neq 0$. The last statement is the 
	consequence of the fact that $B$ is full-rank due to Assumption~\ref{ass:gradI}.
	Therefore, the constraint matrix $C$ is symmetric positive-definite.
\end{proof}

The following theorem gives the reduced-order equations for $\vc q(t)$ such that the quantities $I_k$ are conserved
along the trajectory of the ansatz.

\begin{thm}\label{thm:qdot_const}
	Let assumptions~\ref{ass:ansatz_map} and~\ref{ass:gradI} hold. If a solution to the
	the constrained optimization problem~\eqref{eq:constJ} exists, it must satisfy
	\begin{equation}
	\dot{\vc q} = M^{-1}(\vc q)\left[ \vc f(\vc q) -\sum_{ k = 1 }^{ m } \lambda_k \nabla I_k(\vc q)\right],
	\label{eq:qdot_const}
	\end{equation}
	where $\pmb \lambda=(\lambda_1,\lambda_2,\cdots,\lambda_m)^\top$ is the unique solution to the linear equation
	\begin{equation}
	C \pmb \lambda = \vc b,
	\label{eq:linsys}
	\end{equation}
	with the constrained matrix $C$ defined in Lemma \ref{lem:ASPD} and
	the components of the vector $\vc b = (b_1,b_2,\cdots,b_m)^\top\in\R^m$ defined as
	\begin{equation*}
	b_i = \langle \grad I_i ,M^{-1} \vc f\rangle.
	\end{equation*}
	The metric tensor $M$ and the vector field $\vc f$ are defined as in Theorem~\ref{thm:qdot_unconst}.
\end{thm}

\begin{proof}
	First, we rewrite the constraints in~\eqref{eq:constJ} in an equivalent form. Taking a time derivative, the constraints can be alternatively written as 
	\begin{equation}
	\frac{\id}{\id t}I_k(\vc q(t)) =\langle \nabla I_k(\vc q),\dot{\vc q}\rangle =0, \quad k = 1,2,...,m.
	\end{equation}
	Introducing the Lagrange multiplier $\pmb \lambda = (\lambda_1,...,\lambda_m)^\top \in \R^m$, we define the augmented cost function,
	\begin{equation}
	\mathcal J_c(\vc q,\dot{\vc q},\pmb\lambda) := \mathcal J(\vc q,\dot{\vc q}) + \sum_{ k = 1 }^{ m } \lambda_k\langle \nabla I_k(\vc q),\dot{\vc q}\rangle.
	\label{eq:AugCostFunc}
	\end{equation}	
	If the minimizer of the constrained optimization problem~\eqref{eq:constJ} exists, the partial derivatives of $\mathcal J_c$ with respect to $q_i$ and $\lambda_k$ must vanish at the minimizer. This yields
	\begin{subequations}
		\begin{equation}
		\nabla_{\dot{\vc q}} \mathcal J+  \sum_{ k = 1 }^{ m } \lambda_k \grad I_k =0,
		\label{eq:augLag_grad1n}
		\end{equation}
		\begin{equation}
		\langle \nabla I_1(\vc q),\dot{\vc q}\rangle = \langle \nabla I_2(\vc q),\dot{\vc q}\rangle= ... = \langle \nabla I_m(\vc q),\dot{\vc q}\rangle = 0.
		\label{eq:augLag_grad2n}
		\end{equation}
	\end{subequations}	
	Since $\nabla_{\dot{\vc q}}\mathcal J= M(\vc q) \dot{\vc q} - \vc f(\vc q)$, equation~\eqref{eq:augLag_grad1n} yields
	\begin{equation}
	\dot{\vc q} = M^{-1}(\vc q)\left[ \vc f(\vc q) -\sum_{ k = 1 }^{ m } \lambda_k \nabla I_k(\vc q)\right].
	\label{eq:proof_qdot}
	\end{equation}
	Substituting this expression into~\eqref{eq:augLag_grad2n}, yields $m$ equations
	\begin{equation}
	\sum_{ k = 1 }^{ m } \lambda_k \langle \grad I_i ,M^{-1} \grad I_k\rangle = \langle \nabla I_i ,M^{-1}\vc f\rangle, \quad i = 1,2,...,m.
	\label{eq:augLag_gradn}
	\end{equation}
	Equation~\eqref{eq:augLag_gradn} can be written as the linear system
	$C\pmb \lambda  = \vc b,$
	where $C$ is the constraint matrix~\eqref{eq:constMat}.
	Lemma \ref{lem:ASPD} guarantees that $C$ is invertible and therefore there exists a unique solution $\pmb \lambda$ to the linear system~\eqref{eq:linsys}. 
	Thus, $\dot{\vc q}$ must satisfy equation~\eqref{eq:proof_qdot} with the Lagrange multiplies $\lambda_k$ solving the linear system~\eqref{eq:linsys}.
	This completes the proof.
\end{proof}

We note that unlike the unconstrained problem (Theorem~\ref{thm:qdot_unconst}), the constrained optimization problem~\eqref{eq:constJ}
is not guaranteed to have a solution. For a minimizer to exist, the level sets $\{\vc q\in \Omega : I_k(\vc q) = I_{k,0}\}$, with $k=1,2,\cdots, m$, must have a nonempty intersection. Otherwise, the problem is over-constrained and a minimizer would not exist. 
Apart from these degenerate situations, the constrained optimization problem has a solution and the ansatz evolves according to the
reduced-order equations~\eqref{eq:qdot_const}. In our numerical experiments, presented in Section~\ref{sec:examples}, we never encountered a degenerate case where the solution does not exist.

\section{Numerical examples}
\label{sec:examples}
In this section, we present three numerical examples. 
We begin with a proof-of-concept example where the ansatz captures the exact solution to a linear advection-diffusion equation.
The other two examples deal with nonlinear PDEs, where the ansatz is not an exact solution, but the reduced-order equations nonetheless capture the important features of the system dynamics.

\subsection{Advection-diffusion equation}
As a proof-of-concept example we first consider a linear advection-diffusion equation and an ansatz solution
which is an exact solution of the PDE for an appropriate choice of the time-dependent ansatz parameters.
We show that the reduced-order equations (Theorem~\ref{thm:qdot_unconst}) reproduce this exact solution.

Consider the linear advection-diffusion equation,
\begin{equation}
\dfrac{ \partial u }{ \partial t } =  -c \dfrac{ \partial u }{ \partial x } + \nu \dfrac{ \partial^2 u }{ \partial x^2 }, \quad 
u(x,0) = A_0 \sin \bigg( \dfrac{x}{L_0} \bigg),
\label{eq:advdiff}
\end{equation}
where $c\in\R$ and $A_0, L_0$ and $\nu$ are positive constants.
Equation~\eqref{eq:advdiff} admits the exact solution
\begin{equation}
u(x,t) = A_0 \exp\bigg[ -\frac{ \nu }{ L_0^2 }t \bigg] \sin \bigg(  \frac{ x-ct }{ L_0 } \bigg),
\label{eq:AD_exact}
\end{equation}
which is a traveling sine wave with a decaying amplitude.
We define the ansatz
\begin{equation}
\hat{u}(x, \vc q(t)) = A(t)\sin \bigg(  \frac{ x }{ L(t) } + \phi(t) \bigg),
\label{eq:AD_ansatz}
\end{equation}
with the time-dependent parameters $\vc q(t) = (A(t), L(t), \phi(t))^\top$. Although the ansatz is linear in the amplitude $A$, 
it is a nonlinear function of the parameters $L$ and $\phi$. 

We choose the initial parameter values $A(0) =A_0$, $L(0)=L_0$, and 
$\phi(0) = 0$, so that the ansatz coincides with the initial condition of the PDE, i.e., $u(x,0) = \hat u(x,\vc q(0))$.
Then the ansatz is  an exact solution of \eqref{eq:advdiff} if the parameter values evolve according to
\begin{equation}
A(t) = A_0 \exp\bigg[ -\frac{ \nu }{ L_0^2 }t \bigg], \quad L(t) = L_0, \quad \phi(t) = - \dfrac{ct}{L_0}.
\label{eq:ad_exact_params}
\end{equation}

Now we show that the reduced-order equations~\eqref{eq:qdot_unconst} exactly reproduce 
the parameter evolution~\eqref{eq:ad_exact_params}.
The cost function~\eqref{eq:costfunctional} for the linear advection-diffusion equation reads
\begin{equation}
\mathcal{J} (\hat{u}) = \frac{1}{2} \int_{0}^{2\pi L_0} \big| \hat u_t + c \hat u_x - \nu  \hat u_{xx} \big|^2 \ \id x,
\label{eq:advdiff_cost}
\end{equation}
where the integral is taken over one period of the initial condition.
The appropriate function space for this problem is the Hilbert space of periodic square integrable functions $L^2_{per}(0,2\pi L_0)$.
In this section, we do not enforce any conserved quantities.
After a straightforward calculation, the reduced-order equations~\eqref{eq:qdot_unconst} read
\begin{equation}
\dot{A} = -\dfrac{\nu}{L^2}A,  \quad \dot{L} = 0, \quad \dot{\phi} = -\dfrac{c}{L},
\end{equation}
which have the exact solution~\eqref{eq:ad_exact_params}. In other words, our reduced-order equations applied to the ansatz~\eqref{eq:AD_ansatz} reproduce the exact solution~\eqref{eq:AD_exact}.

\subsection{Nonlinear Schr\"odinger equation}
In this section, we derive reduced-order equations approximating the solutions to the nonlinear Schr\"odinger equation (NLSE).
NLSE is a perturbative model for optical waves~\cite{berge1994,agrawal2013} and surface water wave~\cite{zakharov68,dysthe08}.
NLSE has been intensely studied since, through modulational instability~\cite{benjamin67}, it can reproduce the self-focusing of optical and water waves which leads to the formation of waves of extreme amplitude, often referred to as rogue waves~\cite{solli2007,chabchoub2011}.

The nonlinear Schr\"odinger equation for unidirectional deep water waves is given by
\begin{equation}
\frac{\partial \tilde{u}}{ \partial \tilde{t} } = - i\frac{ \omega_0 }{ 8k_0^2 } \frac{\partial^2 \tilde{u} }{ \partial \tilde{x}^2 } - i\frac{ \omega_0 k_0^2}{ 2 } |\tilde{u}|^2\tilde{u},
\label{eq:FullNLS}
\end{equation}
where $\tilde{u}(\tilde x,\tilde t)$ is the complex wave envelope. The prescribed constants $\omega_0$ and $k_0$ denote the frequency and wave number of the carrier wave, respectively. The wave surface elevation is then given by $\eta(\tilde x,\tilde t) = \mbox{Re}\left[ \tilde u(\tilde x,\tilde t) \exp(k_0\tilde x-\omega_0\tilde t)\right]$.
Introducing the non-dimensional variables $x = 2\sqrt{2} k_0\tilde{x}$,  $t = - \omega_0 \tilde{t}$, and  $u =  ( k_0/\sqrt{2} ) \tilde{u}$, equation~\eqref{eq:FullNLS} becomes 
\begin{equation}
\frac{\partial u}{ \partial t } =   i\frac{\partial^2 u }{ \partial x^2 } + i|u|^2u.
\label{eq:NLS}
\end{equation}
In the following, we work with this non-dimensionalized NLSE.
\begin{figure}[t]
	\centering
	\includegraphics[width=.9\textwidth]{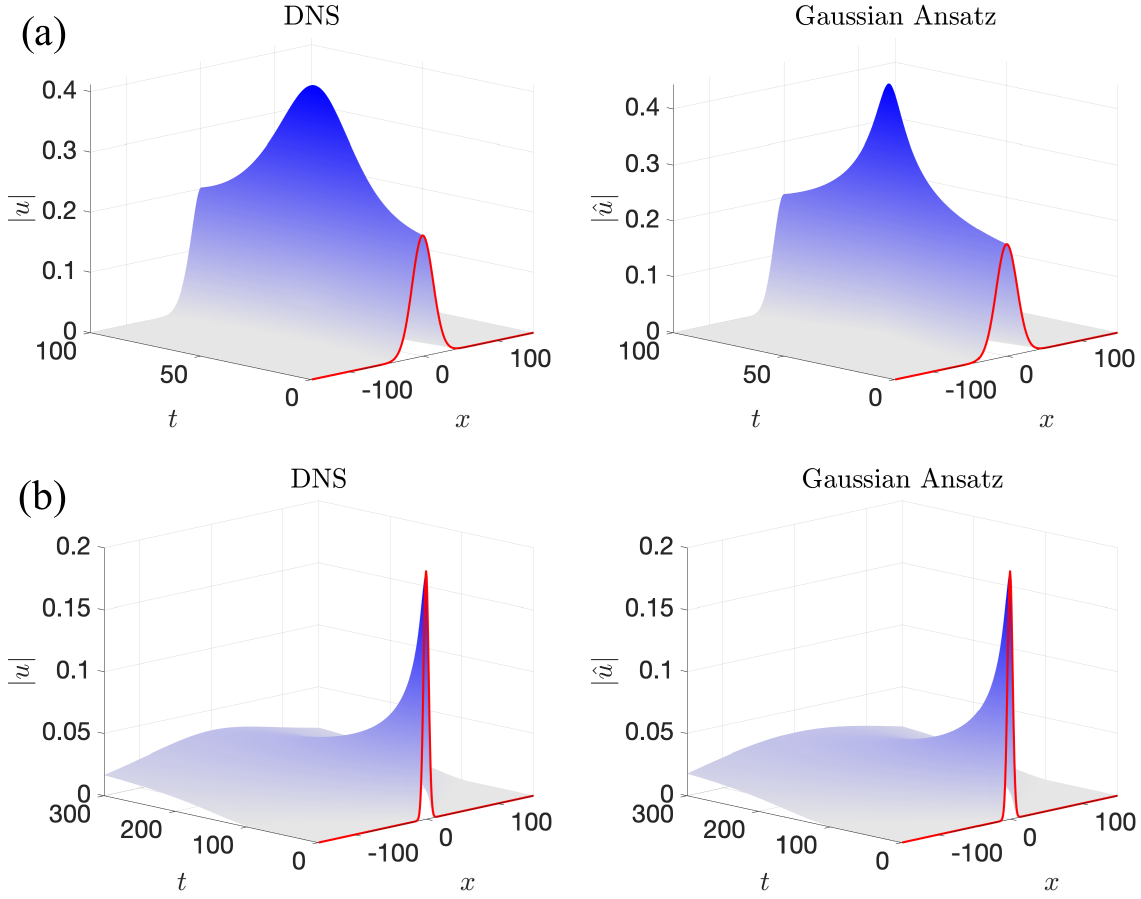}
	\caption{Comparing NLSE solutions using direct numerical simulations (DNS) and the Gaussian ansatz. (a) Focusing wave with initial parameters $A_0 = 0.2$, $L_0 = 20$, $V_0 = -0.05$, $\phi_0 = 0$. (b) Defocusing wave with initial parameters $A_0 = 0.2$, $L_0 = 5$, $V_0 = \phi_0 = 0$.}
	\label{fig:NLS_pcolor}
\end{figure}

It is known that wave groups with certain combinations of initial amplitude and length scale evolve under NLSE to grow in amplitude while their width shrinks. This phenomena is referred to as focusing of wave groups. In contrast, defocusing wave groups exhibit decaying amplitude and widening width (see figure~\ref{fig:NLS_pcolor}).

This observation has spurred a plethora of reduced-order methods for modeling and prediction of rogue waves~\cite{adcock09,Adcock12,ruban2015,ruban2015b,cousins16,PRE2016,farazmand2017}. 
These methods invariably assume a prescribed spatial shape 
for a wave group with time-dependent parameters, such as length scale, amplitude, velocity, and phase. The evolution of these parameters
is then determined by ad hoc methods. Adcock et al.~\cite{adcock09,Adcock12}, for instance, use conserved quantities of NLSE to determine the temporal evolution. Ruban~\cite{ruban2015,ruban2015b} leverages the variational formulation of NLSE to derive a set of ODEs for the parameters. Cousins and Sapsis~\cite{cousins15} take an additional time derivative of NLSE and project the resulting equation
onto their ansatz.

Here, we use the rigorous reduced-order equations~\eqref{eq:qdot_const} to evolve the time-dependent parameters
and show that they correctly capture focusing and defocusing of the wave groups.
Following earlier work~\cite{PerezGarcia1996,ruban2015,ruban2015b}, we use the Gaussian ansatz,
\begin{equation}
\hat{u}(x, \vc q (t)) = A(t) \exp \bigg[ -\frac{x^2}{L^2(t)} + i \frac{ x^2 V(t) }{L(t)} + i \phi(t) \bigg],
\label{eq:NLS_Ansatz}
\end{equation}
with the time-dependent parameters $\vc q(t) = (A ,L, V, \phi)$, where 
$A$ determines the wave amplitude, $L$ is a length scale controlling how quickly the wave group 
decays away from its center, $V$ is the wave velocity, and $\phi$ is the wave phase. 
The motivation for choosing this ansatz is that, for a special choice of the parameters $\vc q(t)$, it is an exact solution to the linear part of NLSE~\cite{ruban2015b}. However, ansatz~\eqref{eq:NLS_Ansatz} is not an exact solution of NLSE except for the trivial case $A(t) = 0$. 
We note that $\hat{u}$ is complex-valued, but its parameters $\vc q = (A ,L, V, \phi)$ are real-valued with the additional restrictions that $A(t) >0$ and $L(t)>0$.

The cost function~\eqref{eq:costfunctional} for NLSE reads
\begin{equation}
\mathcal{J}( \vc q, \dot{ \vc q } ) = \frac{1}{2} \int_{\R} \big|\hat{u}_t - i\hat{u}_{xx} - i |\hat{u}|^2\hat u \big|^2 \ \id x,
\end{equation}
where $|\cdot|$ denotes the modulus of a complex number. The Hilbert space $H$ here is the square integrable complex functions
over $\R$. We note that, for simplicity, we stated our results in Section~\ref{sec:varmodeling} for real-valued functions. Their generalization to complex-valued functions is straightforward.

NLSE has several conserved quantities~\cite{adcock09,Adcock12} which can be enforced using the method described in 
Section~\ref{sec:conserved_quants} (Theorem~\ref{thm:qdot_const}). Here, we only enforce the conservation of the two most relevant quantities,
\begin{equation}
I_1(u) = \int_{\R} |u|^2 \ \id x,\quad 
I_2(u) = \int_{\R} |u_x|^2 \ \id x - \frac{1}{4}\int_{\R} |u|^4 \ \id x,
\end{equation} 
which are the total mass and the total energy, respectively.
To obtain the reduced-order equations~\eqref{eq:qdot_const}, 
we solve the constrained optimization problem~\eqref{eq:constJ} with $m=2$ and
$$I_{1,0} = \sqrt{\frac{\pi }{2}} A_0^2 L_0 , \quad I_{2,0}= \frac{\sqrt{\pi } A_0^2 \left(2 \sqrt{2} \left(L_0^2 V_0^2+1\right)-A_0^2 L_0^2\right)}{8 L_0} ,$$
where $(A_0,L_0,V_0,\phi_0)$ are the initial values of the parameters of the ansatz~\eqref{eq:NLS_Ansatz}.
 The resulting reduced-order ODEs  read
	\begin{equation}
	\dot{ A } = -\frac{2 A V}{L}, \quad \dot{ L } = 4 V, \quad  \dot{ V } =  \frac{4}{L^3}-\frac{A^2}{\sqrt{2} L}, \quad \dot{ \phi } =  \frac{5 A^2}{4 \sqrt{2}}-\frac{2}{L^2}.
	\label{eq:NLS_RONS_ODEs}
	\end{equation}
\begin{figure}[t]
	\centering
	\includegraphics[width=.49\linewidth]{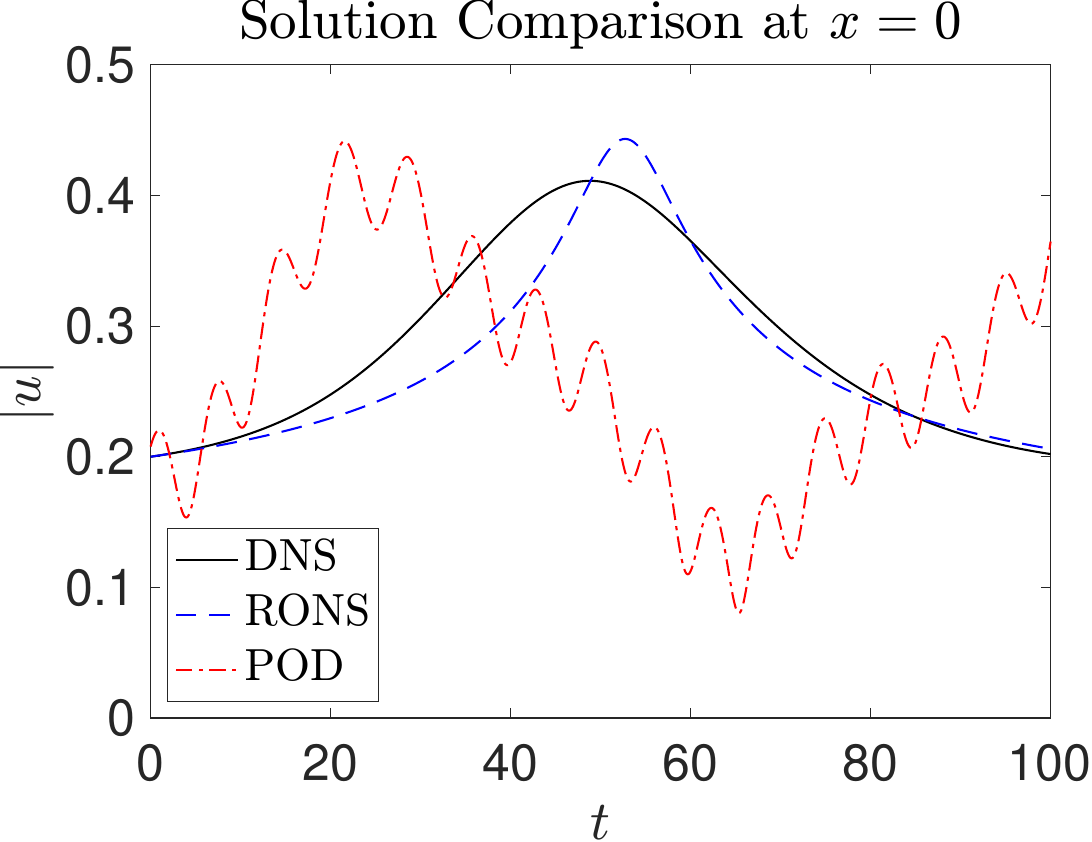}
	\includegraphics[width=.49\linewidth]{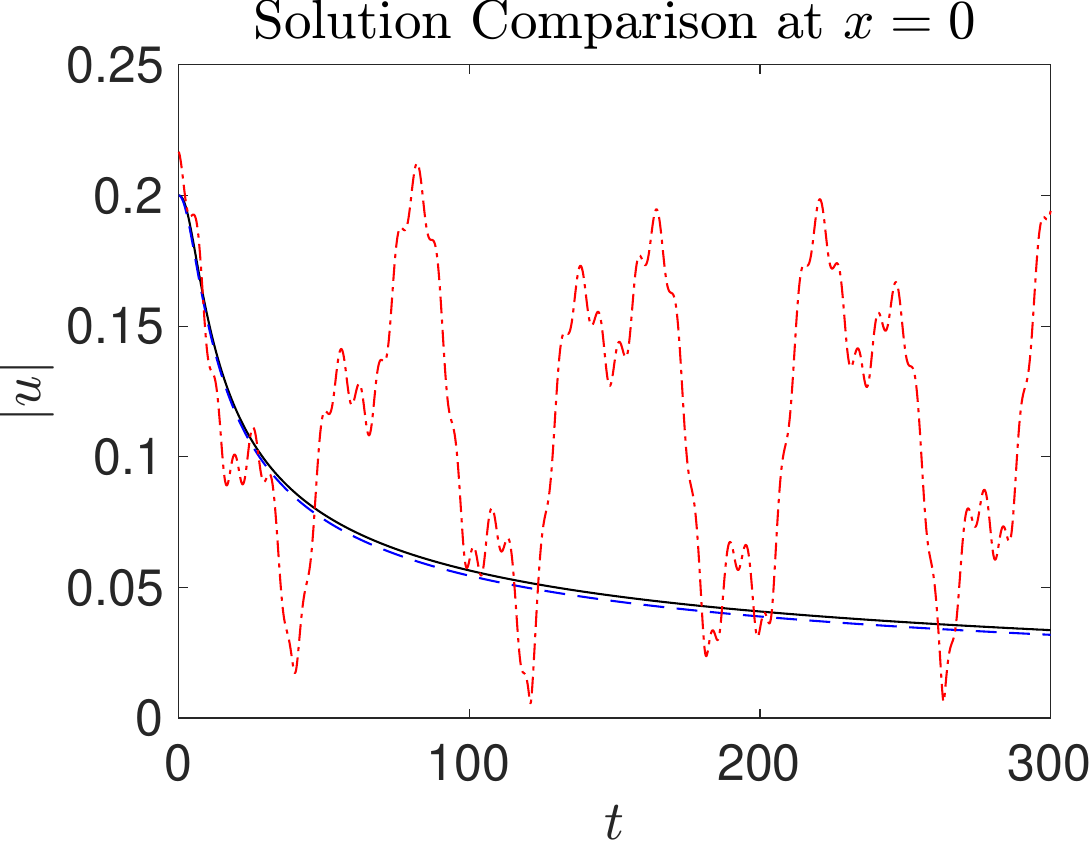}
	\caption{  Evolution of NLSE solutions evaluated at $x=0$ using direct numerical simulations (DNS), RONS with the Gaussian ansatz, and Galerkin projection onto four POD modes. Left:
		Focusing wave with initial parameters $A_0 = 0.2$, $L_0 = 20$, $V_0 = -0.05$, $\phi_0 = 0$. Right: Defocusing wave with initial parameters $A_0 = 0.2$, $L_0 = 5$, $V_0 = \phi_0 = 0$.}
	\label{fig:NLS_center}
\end{figure}

For comparison, we use direct numerical simulations (DNS) of the NLSE as the ground truth.
As in Ref.~\cite{cousins16}, we solve NLSE using a Fourier pseuedospectral method in space and a fourth-order Runge–Kutta exponential time differencing scheme \cite{cox02}.
Boundary conditions are assumed to be periodic for ease of implementation and should not affect the results so long as the domain is chosen to be large enough so that $u$ is small near the boundaries. 
We use the domain size $256\sqrt{2}\pi$ ($64$ wave periods) which we have found large enough to ensure $u$ is small near the boundary for the duration of the computation.
For all runs we use $2^{10}$ Fourier modes and a time step of 0.025. Using more Fourier modes, smaller time steps, or larger domain size did not significantly affect the results of our simulations. 
For example, doubling the number of Fourier modes and halving the time step caused a relative change of $|u|$ on the order of $10^{-7}$  in the $L^2$-norm.
The DNS results are initialized with $u_0(x) = \hat u(x,\vc q(0))$, so that the DNS solution and the ansatz solution coincide at the initial time.

In figure~\ref{fig:NLS_pcolor}, we plot the evolution of the DNS results $|u|$ and the ansatz $|\hat u|$ for two sets of initial parameters. One set of parameters leads to a focusing wave (growing amplitude) and the other leads to a defocusing wave (decaying amplitude).
In both cases, the ansatz correctly predicts the qualitative focusing or defocusing behavior of the wave group. We have repeated these comparisons for a range of parameter values $\vc q_0$ and in every case, the ansatz correctly predicts the focusing or defocusing behavior of the wave group.

In figure \ref{fig:NLS_center}, we plot the evolution of the amplitudes at $x=0$ for the same set of parameter values as in figure~\ref{fig:NLS_pcolor}. 
For the case of the defocusing wave, the reduced-order solution produced by RONS is in excellent agreement with the DNS results.
For the focusing wave the RONS reduced-order solution provides a reasonable approximation compared to the DNS, roughly capturing both the peak amplitude and the time that it occurs. However, the Gaussian ansatz overestimates the peak height and also produces a peak which is thinner than that of the DNS. This behavior, i.e. excellent agreement for defocusing waves and overshooting for focusing waves, is systematically observed for other parameter values $\vc q_0$ (not shown here).

 To compare RONS results with a standard technique in model reduction, figure~\ref{fig:NLS_center} also shows the results obtained by Galerkin projection with modes found through proper orthogonal decomposition (POD) as described in \cite{shlizerman2011}. 
We use four POD modes as a benchmark which is on par with the four parameters $(A,L,V,\phi)$ involved in the Gaussian RONS ansatz~\eqref{eq:NLS_Ansatz}.
We point out that the POD modes are complex valued, so that the reduced-order equations involve four complex-valued ODEs or equivalently eight real-valued ODEs.
In other words, the POD-reduced equations involve twice as many ODEs as the RONS reduced-order equations~\eqref{eq:NLS_RONS_ODEs}.
Nonetheless, as shown in figure~\ref{fig:NLS_center}, the POD reduced-order solutions provide poor approximations for both focusing and defocusing waves, 
demonstrating the superiority of RONS.
To obtain an accurate POD-reduced model which is capable of accurately approximating the focusing and defocusing waves, we have observed that at least 16 complex-valued POD modes are required.


\subsection{Two-dimensional fluid flow}
For our last example, we consider the flow of a  two-dimensional, incompressible, and Newtonian fluid. 
The vorticity equation for such a fluid is given by
\begin{equation}
\pard{\omega}{t} +  \vc u \cdot \grad \omega =  \nu \Delta\omega,
\label{eq:Vorticity_Eqn}
\end{equation}
where $\nu$ is the kinematic viscosity, $\vc u(\vc x,t)$ denotes the velocity field of the fluid, and $\omega(\vc x,t)$ is the
component of the vorticity orthogonal to the plane of motion. We consider the flow on the unbounded two-dimensional domain, $\vc x = (x, y)\in\R^2$. 
A two-dimensional incompressible fluid admits a stream function $\psi(x, y,t)$ which satisfies 
\begin{equation}
\vc u(\vc x, t) = ( \psi_y, -\psi_x )^\top,\quad \omega = -\Delta \psi.
\end{equation}

Therefore, prescribing the stream function $\psi$ is sufficient to determine the fluid velocity $\vc u$ and the vorticity $\omega$.  
As a result, we define the ansatz in terms of the stream function,
\begin{equation}
\hat \psi(\vc x, \vc q(t)) = \sum_{i = 1}^{N} A_i(t) \exp \bigg[ -\frac{ | \vc x - \vc x_{i}(t) |^2}{L_i^2(t)} \bigg],
\label{eq:psi_anzats}
\end{equation} 
which is the superposition of $N$ axisymmetric vortices,
where $A_i(t)$ denotes the amplitude of the $i$-th vortex, $L_i(t)$ is its length scale, and $\vc x_{i}(t) = (x_{i}(t), y_{i}(t))$ is the vortex center. 
Therefore, each term of the ansatz has four time-dependent parameters, $(A_i,L_i,x_i,y_i)$, resulting in a total of $n=4N$ parameters
for $N$ vortices.
We use this ansatz because it belongs to the class of smooth vortex methods~\cite{chorin_1973,koumoutsakos_2000} and, more importantly, vortices with the Gaussian stream function~\eqref{eq:psi_anzats} have been observed in laboratory experiments~\cite{trieling1997}. 
\begin{figure}[t]
	\centering
	\includegraphics[width=.7\linewidth]{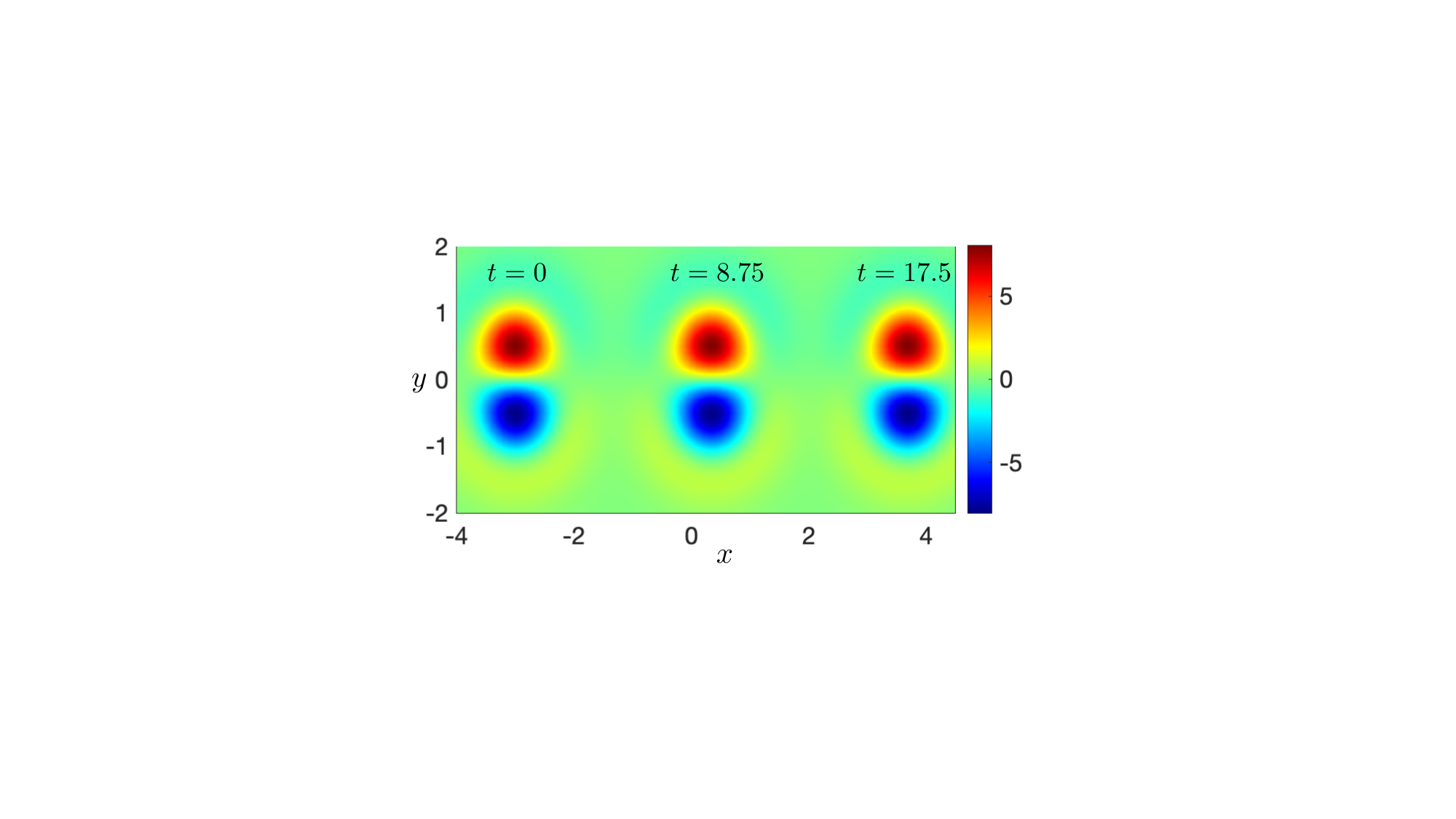}
	\caption{Reduced dynamics of a vortex dipole. The initial parameter so the ansatz are $A_1(0) = 1$, $A_2(0) = -1$, $L_1(0) = L_1(0) = 0.75$, $\vc x_1(0)  = (-3,0.5)$, $\vc x_2(0)  = (-3,-0.5)$. See the supplementary material for a movie.}
	\label{fig:TwoVortex_opposite}
\end{figure}

After defining the stream function, the corresponding fluid velocity $\hat{\vc u}$ and vorticity $\hat\omega$ are computed to define the cost functional,
\begin{equation}
\mathcal{J} ( \vc q, \dot{ \vc q }) =  \frac{1}{2} \int_{\R^2} | \hat \omega_t + \hat{\vc u} \cdot \grad \hat \omega-\nu\Delta\hat \omega|^2 \ \id \vc x.
\label{eq:fluid_costfunctional}
\end{equation}
For brevity, here we only study the inviscid case, $\nu=0$, where the vortex dynamics has a rich mathematical structure.
In particular, the inviscid vorticity equation admits a class of exact weak solutions called \emph{point vorticies}~\cite{arefVortex}.
We show that the reduced equations governing the parameters of the ansatz~\eqref{eq:psi_anzats} reproduce the expected
point vortex dynamics, although the ansatz is not an exact solution of the vorticity equation.

Point vortex solutions of the inviscid vorticity equation consist of the superposition of Dirac delta functions,
\begin{equation}
\omega (\vc x,t) = \sum_{i=1}^N \Gamma_i \delta (\vc x -\vc x_i(t)),
\label{eq:w_pv}
\end{equation}
where $\Gamma_i$ are constant vortex amplitudes and $\vc x_i(t) = (x_i(t),y_i(t))$ are the time-dependent vortex centers.
Using the Green's function for the Poisson equation $\Delta \psi = -\omega$, the stream function for point vortex solution~\eqref{eq:w_pv}
can be obtained as
\begin{equation}
\psi(\vc x,t) = -\frac{1}{2\pi}\sum_{i=1}^N \Gamma_i \log |\vc x-\vc x_i(t)|.
\label{eq:psi_pv}
\end{equation}
For the point vortex to be a weak solution of the vorticity equation, the vortex centers $\vc x_i(t)$ must 
satisfy the set of ordinary differential equations~\cite{newton_Nvortex},
\begin{equation}
\Gamma_i\dot x_i = \pard{\mathcal H}{y_i},\quad \Gamma_i\dot y_i= - \pard{\mathcal H}{x_i},
\end{equation}
where $\mathcal H = -\sum_{i\neq j} \Gamma_i\Gamma_j\log |\vc x_i-\vc x_j|/4\pi$ is the Hamiltonian corresponding to the 
$N$ vortex problem.

Dynamics of point vortices depends on the number of vortices $N$ and their initial configuration. 
We refer the reader to Refs~\cite{arefVortex,newton_Nvortex} for a complete accounting.
Here, we consider three distinct configurations,
\begin{enumerate}
	\item Vortex dipole: Two vortices of equal strengths but with opposite signs,
	\item Vortex pair: Two vortices of equal strengths and signs,
	\item Leapfrogging: Two vortex pairs with opposite signs.
\end{enumerate}
For each configuration, we show that the reduced-order equations~\eqref{eq:qdot_const} applied to the ansatz~\eqref{eq:psi_anzats}
reproduce the expected vortex dynamics predicted by the point vortex solution. We emphasize that, unlike the point vortex~\eqref{eq:psi_pv}, our ansatz has a 
smooth profile and is not an exact solution of the vorticity equation. 
In addition, we allow the amplitude $A_i$ and the length scale $L_i$ of each vortex in the ansatz to vary with time.

The inviscid vorticity equation admits several conserved quantities~\cite{newton_Nvortex}. 
Here, we enforce two of these conserved quantities, 
\begin{equation}
I_1(\vc u) = \frac{1}{2} \int_{\R^2}| \vc u |^2 \ \id \vc x,\quad 
I_2(\vc u)=\frac{1}{2} \int_{\R^2}  | \nabla\times \vc u |^2 \ \id \vc x = \frac{1}{2} \int_{\R^2}  |  \omega |^2 \ \id \vc x,
\end{equation}
which denote the kinetic energy and enstrophy, respectively.
To derive the reduced-order equations, we minimize the cost function~\eqref{eq:fluid_costfunctional} with the constraints $I_1 = I_{1,0}$
and $I_2 = I_{2,0}$ (see Theorem~\ref{thm:qdot_const}).

First, we present the results for a vortex dipole. Point vortex dynamics predict that a vortex dipole, consisting of two vortices of opposite sign, travels on a straight line without changing shape.
Figure \ref{fig:TwoVortex_opposite} shows a similar scenario for our smooth Gaussian ansatz, where
two vortices of opposite sign and equal amplitude are placed symmetrically along the $x$-axis. 
The vortices are evolved according to our reduced-order equation~\eqref{eq:qdot_const}.
The length-scales $L_i$ and amplitudes $A_i$ remain constant during the evolution. The vortex centers $\vc x_i$, however,
translate horizontally at a constant speed, reproducing the expected dynamics that point vortices exhibit.
\begin{figure}[t]
	\centering
	\includegraphics[width=\textwidth]{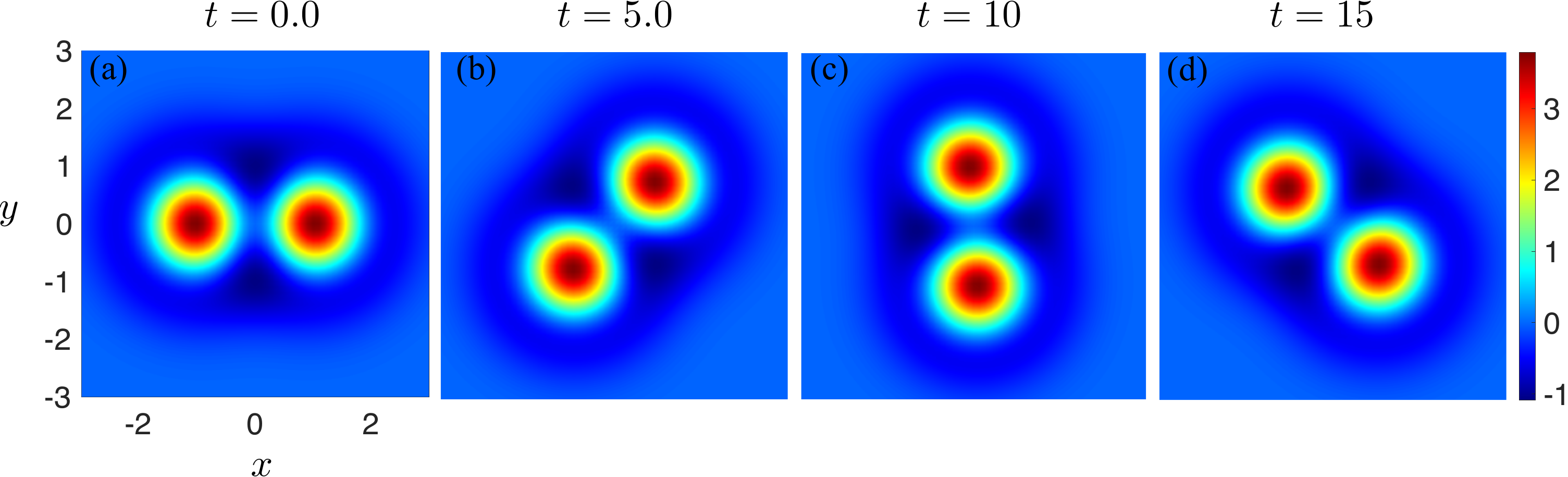}
	\caption{Reduced dynamics of a vortex pair. The initial parameter so the ansatz are $A_1(0) = A_2(0) = 1$, $L_1(0) = L_2(0) = 1$, $\vc x_1(0)  = (-1,0)$, $\vc x_2(0)  = (1,0)$. See the supplementary material for a movie.}
	\label{fig:TwoVortex_same}
\end{figure}

Next, we consider the case of a vortex pair, where the vortices have the same sign and amplitudes.
Point vortex dynamics, in this case, predicts that the vortices rotate around their midpoint at a constant angular velocity. 
Figure \ref{fig:TwoVortex_same} shows a similar vortex pair created by the Gaussian ansatz~\eqref{eq:psi_anzats}.
The reduced-order dynamics shows that the system is in a relative equilibrium where the two vortices continually rotate counterclockwise around the origin $\vc x=0$ at a constant angular velocity. Again, the reduced-order dynamics are in excellent agreement with the expected point vortex behavior.
\begin{figure}[t!]
	\centering
	\includegraphics[width=\textwidth]{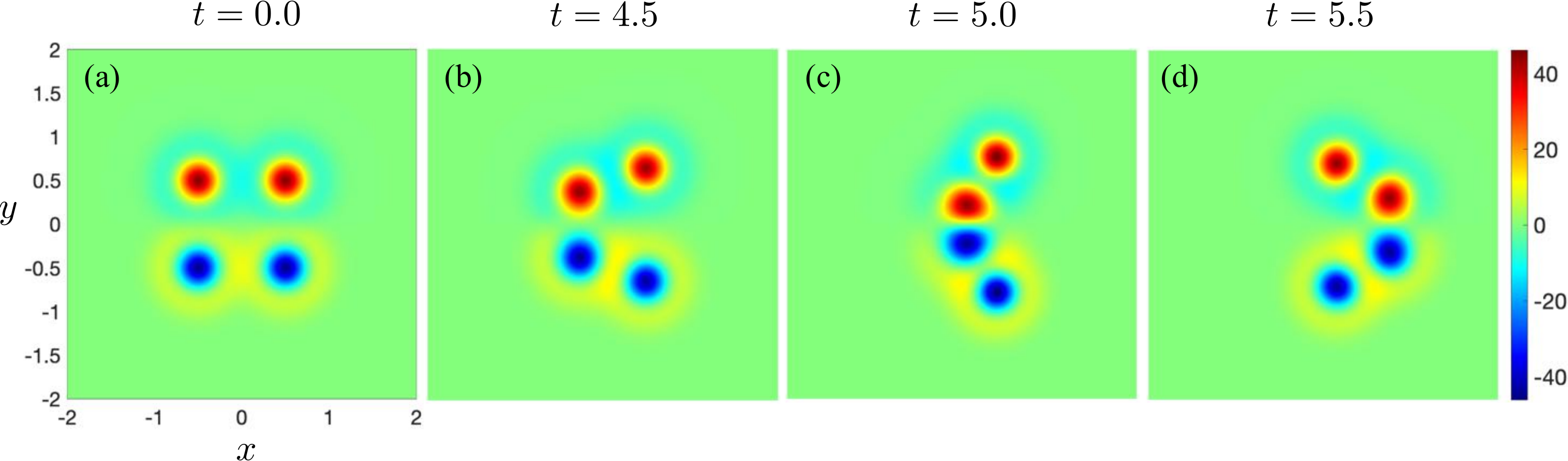}
	\caption{Reduced dynamics of leapfrogging vortices. The initial parameters values are $A_1(0) = A_3(0) = 1$, $A_2(0) = A_4(0) = -1$, $L_i(0) = 0.3$, $\vc x_1(0) = (0.5, 0.5)$ , $\vc x_2(0) = (0.5, -0.5)$  , $\vc x_3(0) = (-0.5, 0.5)$, $\vc x_4(0) = (-0.5, -0.5)$. See the supplementary material for a movie.}
	\label{fig:FourVortex_Leap}
\end{figure}

Finally, we consider the leapfrogging configuration, where four vortices are placed at the corners of a square (see figure~\ref{fig:FourVortex_Leap}).
The vortices on the top row have positive vorticity whereas the vortices on the bottom row have negative vorticity.
The point vortex dynamics predicts that the vortices to the left accelerate and zoom in between the vortices in the front. 
This motion repeats indefinitely, creating the so-called leapfrogging dynamics.
Figure \ref{fig:FourVortex_Leap} shows the reduced dynamics of the Gaussian ansatz~\eqref{eq:psi_anzats}.
One pair of vortices is placed at $(\pm0.5, 0,5)$, and another pair of vortices with opposite sign is placed at $(\pm 0.5, - 0,5)$.
The RONS reduced-order equations once again reproduce the leapfrogging dynamics predicted by point vortices.  The computational time to numerically integrate the RONS equations for the vortex dipole, vortex pair, and leapfrogging configurations was $0.0087$, $0.0304$, and $4.2188$ seconds, respectively, on a 2019 Macbook Pro with a 1.7 GHz Quad-Core Intel Core i7 processor.

\begin{figure}
	\centering		
	\includegraphics[width=\linewidth]{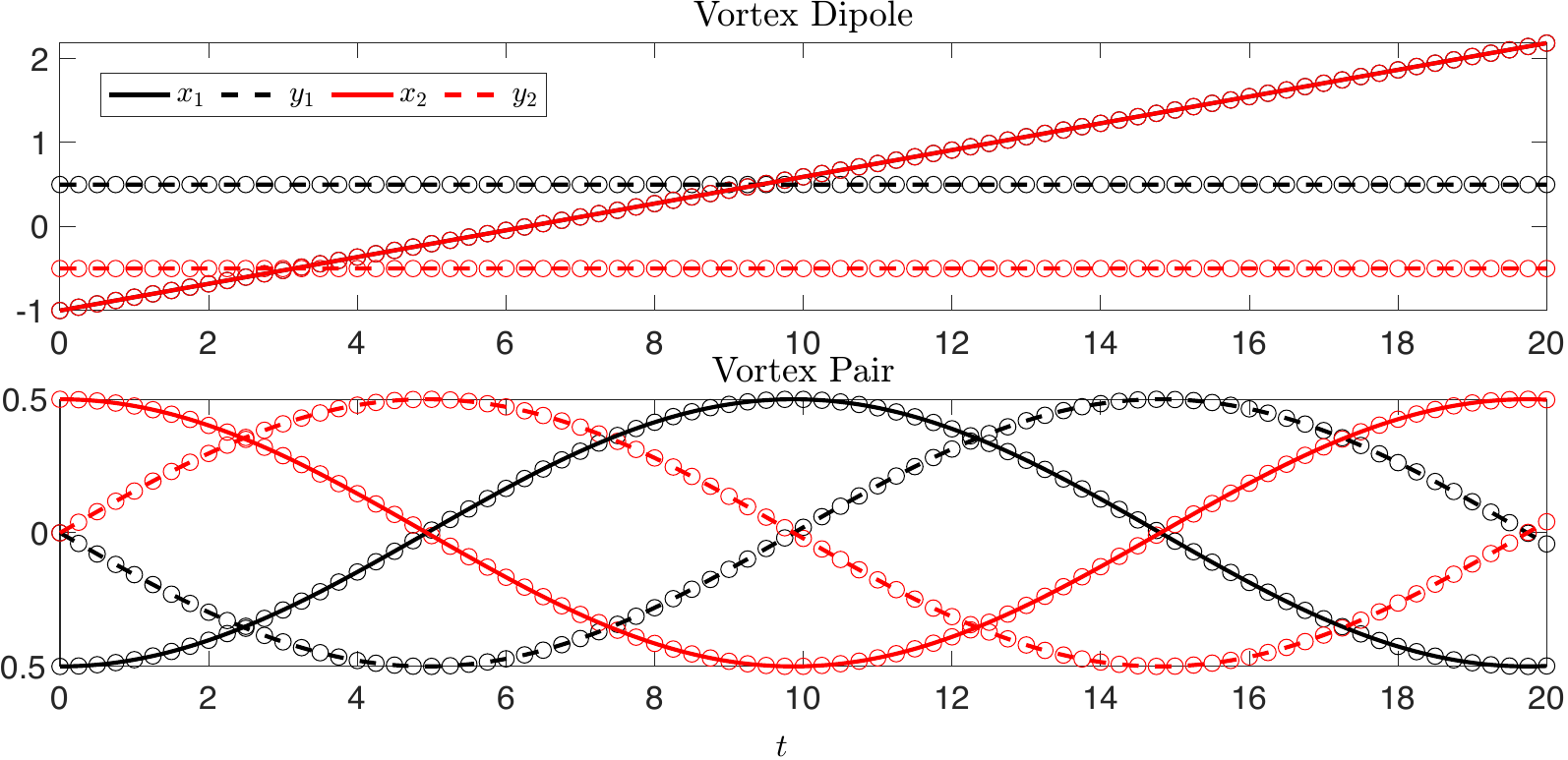}
	\caption{  Comparing RONS against point vortex dynamcis for a vortex dipole (top) and a vortex pair (bottom). RONS trajectories are solid and dashed lines and point vortex trajectories are marked by circles. The parameter values used for the vortex dipole are $\Gamma_1 = 1$, $\Gamma_2  =  -1$, $\eps(0) = 0.05$, $\vc x_1(0) =  (-1, 0.5)$, $\vc x_2(0) =  (-1, -0.5)$, and initial parameter values for the vortex pair are  $\Gamma_1 = \Gamma_2  =  1$,  $\eps(0) = 0.05$, $\vc x_1(0) =  (-0.5, 0)$,  $\vc x_2(0) =  (0.5, 0)$. Note that, in the vortex dipole plot, the time series of $x_1$ and $x_2$ overlap. }
	\label{fig:Pointvortex_Comparison}
\end{figure}

 We chose the Gaussian stream function~\eqref{eq:psi_anzats} as the ansatz since it has been observed in experiments~\cite{trieling1997}; then we showed that the RONS reduced-order equations reproduce the expected dynamics. However, this ansatz is not convenient for a quantitative comparison with solutions produced by point vortices. 
To make quantitative comparisons with point vortex dynamics, we introduce the following ansatz for the vorticity,
\begin{equation}
	\hat \omega(\vc x, \vc q) =  \sum_{i=1}^N \frac{\Gamma_i}{2\pi \eps^2}\exp \left[ - \frac{ |\vc x-\vc x_i|^2 }{ 2\eps^2 } \right],
	\label{w_RONS_pv}
\end{equation}
which is a sum of Gaussians with the amplitudes $A_i = \Gamma_i / (2\pi \eps^2)$ and lengthscales $L_i=\epsilon\sqrt{2}$. Note that the ansatz~\eqref{w_RONS_pv} converges to  point vortices~\eqref{eq:w_pv} as $\eps$ tends to zero.
After obtaining the corresponding fluid velocity $\hat{\vc u}$ using the Biot-Savart law~\cite{newton_Nvortex}, we derive the RONS equations by minimize the cost functional~\eqref{eq:fluid_costfunctional} while enforcing conservation of kinetic energy and enstrophy. Figure \ref{fig:Pointvortex_Comparison} shows the vortex center trajectories with $\epsilon(0)=0.05$ for the vortex dipole and the vortex pair. It shows that the RONS results are in excellent agreement with point vortex dynamics. 
We point out that, to derive the RONS reduced-order equations, we allow $\epsilon(t)$ and $\vc x_i(t)$ to be time-dependent functions. Nonetheless, RONS predicts that $\epsilon$ remains constant ($\dot\epsilon=0$).

In closing, we point out that our reduced-order equations 
can be similarly applied to the viscous vorticity equation ($\nu>0$) with no difficulty.
In this case, the equation is dissipative and therefore there are no conserved quantities to be enforced.
However, direct numerical simulations indicate that initially axisymmetric vortices evolve under the viscous vorticity equation
into an approximately elliptical shape. As a result, the axisymmetric ansatz~\eqref{eq:psi_anzats} may no longer be appropriate, 
and a more general non-axisymmetric Gaussian ansatz might be needed. This viscous case will be explored in future work.

\section{Conclusions}
\label{sec:conclusion}
We proposed RONS as a general framework for evolving time-dependent nonlinear reduced-order solutions of PDEs which is applicable to a broad class of problems. The reduced-order solution has a prescribed shape in space and depends nonlinearly on a set of time-dependent parameters.
Our reduced-order equations evolve the parameters such that the instantaneous error between the ansatz dynamics and the full PDE dynamics is minimized. Any number of conserved quantities can readily be enforced in our framework without requiring a Lagrangian or Hamiltonian formulation for the PDE.

We presented several numerical examples to demonstrate the effectiveness of this method: a linear advection-diffusion equation, the nonlinear Schr\"odinger equation, and Euler's equation for two-dimensional ideal fluids.
In every case, the reduced-order equations produce approximate solutions which capture all qualitative features of the PDE.

Our method opens a new paradigm in reduced-order modeling, with many  aspects undoubtedly remaining to be explored.
For instance, an upper bound should be derived to estimate the accumulated finite-time error between the reduced-order solution
and the true solutions of the PDE.
Here we focused on approximating strong solutions; further work is needed to derive reduced-order equations for PDEs in the weak 
formulation. Finally, here we resolved an important aspect of reduced-order nonlinear solutions, namely, the evolution of the ansatz variables $\vc q(t)$. 
However, the choice of the shape of the ansatz $\hat u(\vc x,\vc q(t))$ still relies on domain expertise and familiarity with the PDE. A systematic method
for determining an appropriate ansatz manifold remains an open problem.

\appendix
\section{Instability of the finite-time formulation}
\label{sec:EL_unstable}
In Remark~\ref{rem:FTerror}, we advised against using the finite-time error,
\begin{equation}
S[\vc q] = \int_0^T \mathcal J (\vc q(t),\dot{\vc q}(t))\id t,
\label{eq:FTerror}
\end{equation}
since it generally leads to 
unstable reduced-order models. In this section, we demonstrate this with a simple example.
Consider the linear PDE,
\begin{equation}
\pard{u}{t} = \mathcal L u,
\end{equation}
where $\mathcal L$ is a self-adjoint negative-definite operator with orthonormal eigenfunctions $u_i(\vc x)$ and corresponding eigenvalues $-\lambda_i$, where $\lambda_i>0$.

Approximating a solution of the PDE as the linear combination of the eigen functions $u_i$, we consider the linear ansatz,
\begin{equation}
\hat{u} (x, \vc q(t)) = \sum_{i=1}^N q_i(t)u_i(x).
\label{eqn:linear_pde_ansatz}
\end{equation}
The standard Galerkin projection leads to the uncoupled ODEs,
\begin{equation}
\dot{q}_k = -\lambda_k q_k,\quad k=1,2,\cdots,N,
\end{equation}
which admit the exact solution $q_k(t) = q_k(0) \exp(-\lambda_k t)$. Therefore, the Galerkin method correctly predicts that the
solutions decay to zero exponentially fast.

Next we derive the reduced-order equations by extremizing the action~\eqref{eq:FTerror}.
Substituting the linear ansatz~\eqref{eqn:linear_pde_ansatz} in the Lagrangian $\mathcal J$, we obtain
\begin{equation}
\mathcal{J} (\vc q, \dot{ \vc q }) = \frac{1}{2} \int_{\R} | \hat u_t - \mathcal L \hat  u |^2 \ \id \vc x 
= \frac{1}{2} \sum_{i=1}^N \left(\dot q_i^2 +2\lambda_i\dot q_i q_i +\lambda_i^2q_i^2\right).
\end{equation}
We take the first variation of the functional $S[\vc q]$ with respect to perturbations $\delta \vc q(t)$, $t\in[0,T]$.
The variations at the endpoints, $\delta \vc q(0)$ and $\delta \vc q(T)$, need to be prescribed. Since, in general, the true solution 
of the PDE is a priori unknown, the endpoint variations cannot be prescribed by taking the different between the ansatz solution
and the true solution. Therefore, we adopt the usual assumption from classical mechanics that the variations vanish at the end points.
Requiring the first variation of $S$ to vanish
leads to the Euler-Lagrange equations,
\begin{equation}
\frac{\id }{ \id t }\frac{\partial \mathcal{J}}{ \partial \dot{q}_k }-\frac{\partial \mathcal{J}}{ \partial q_k } =
\ddot q_k -\lambda_k^2q_k=0,
\end{equation}
which admit the exact solution,
\begin{equation}
q_k(t) = \left[ \frac12 q_k(0) + \frac{1}{2\lambda_k}\dot q_k(0)\right] e^{\lambda_k t} +  \left[ \frac12 q_k(0) - \frac{1}{2\lambda_k}\dot q_k(0)\right] e^{-\lambda_k t}.
\end{equation}
Except for the special initial conditions $\dot q_k(0)=-\lambda_kq_k(0)$, the solutions grow exponentially fast.
However, we know that the true solutions of the PDE must decay to zero. 
Even if the special initial condition $\dot q_k(0)=-\lambda_kq_k(0)$ is specified,
numerical round-off errors will grow exponentially in time, rendering the solution unstable in practice.
In the case of the nonlinear Schr\"odinger equation,
we have observed (not presented here) a similar instability for reduced-order equations obtained 
from the finite-time error functional~\eqref{eq:FTerror}.


\end{document}